\def\Bbb#1{{\bf #1}}

\def\fnote#1{\footnote}
\def\blacksquare{\hbox{\vrule width 4pt height 4pt depth 0pt}}

\def\cwleftpar#1#2{\leftskip #1 \rightskip #2 plus 1fill}
\def\cwrightpar#1#2{\leftskip #1 plus 1fill \rightskip #2}
\def\cwcenterpar#1#2{\leftskip #1 plus 1fill \rightskip #2 plus 1fill}
\def\cwfullpar#1#2{\leftskip#1\rightskip#2}

\def\cwoutdent#1#2{\llap{\hbox to #1{#2 \hss}}\ignorespaces}
\def\cwparbegin#1#2#3#4#5{
	\ifcase #1 \cwleftpar{#2}{#3}
	\or \cwrightpar{#2}{#3}
	\or \cwcenterpar{#2}{#3}
	\else \cwfullpar{#2}{#3}\fi
	\ifcase #4 \baselineskip = 1.5\baselineskip
	\or \baselineskip = 2\baselineskip
	\or \baselineskip = 3\baselineskip
	\else \baselineskip = 1\baselineskip\fi
	\ifdim #5 > 0in \else \noindent \fi
	\noindent\ignorespaces}
\documentclass{article}
\begin{document}
\advance \vsize by -1\baselineskip
\def\makefootline{
{\vskip \baselineskip \noindent \folio                               \par
}}

\vspace*{2ex}
\noindent {\Huge Linear Transports along Paths in\\[0.4ex] Vector
Bundles}\\[1.5ex]
\noindent  {\Large I. General Theory}

\vspace*{2ex}
\noindent Bozhidar Zakhariev Iliev
\fnote{0}{\noindent $^{\hbox{}}$Permanent address:
Laboratory of Mathematical Modeling in Physics,
Institute for Nuclear Research and \mbox{Nuclear} Energy,
Bulgarian Academy of Sciences,
Boul.\ Tzarigradsko chauss\'ee~72, 1784 Sofia, Bulgaria\\
\indent E-mail address: bozho@inrne.bas.bg\\
\indent URL: http://theo.inrne.bas.bg/$\sim$bozho/}

\vspace*{2ex}

{\bf \noindent Published: Communication JINR, E5-93-239, Dubna, 1993}\\[1ex]
\hphantom{\bf Published: }
http://www.arXiv.org e-Print archive No.~math.DG/0411023\\[2ex]

2003 PACS numbers: 02.40.Ma, 02.40.Vh, 04.90.+e\\[2ex]

\noindent
{\small
The \LaTeXe\ source file of this paper was produced by converting a
ChiWriter 3.16 source file into
ChiWriter 4.0 file and then converting the latter file into a
\LaTeX\ 2.09 source file, which was manually edited for correcting numerous
errors and for improving the appearance of the text.  As a result of this
procedure, some errors in the text may exist.
}\\[2ex]

	\begin{abstract}
The (parallel) linear transports along paths in vector bundles are
axiomatically described. Their general form and certain properties are found.
It is shown that these transports are locally (i.e. along every fixed path)
always Euclidean ones in a senses that there exist frames in which their
matrices are unit. The investigated transports along paths are described in
terms of their local coefficients, as well as in terms of derivations along
paths.
	\end{abstract}\vspace{3ex}

 {\bf 1. INTRODUCTION}

\medskip
The corresponding to linear connections parallel transports are linear isomorphisms between the tensor spaces along the paths they act [1,2]. In the present work we appropriately generalize these transports, preserving the property that they are linear isomorphisms between vector spaces. As a ground for this we take the basic properties of the considered in [3] (parallel) transports along paths generated by derivations of tensor algebras. In some places we also follow the ideas of this reference by transferring mutatis mutandis definitions and results. Part of the material of the present investigation was realized, in the case of the tangent bundle of a differentiable manifold, in [4]. That reference also contains possible applications to physics which, through some changes, have a place also with respect to this work.

In Sect. 2, we define and find the general, invariant and in components, form
of the linear transports along paths in vector bundles. In Sect. 3 is proved
the existence of frames in which these transports look like an usual parallel
transport in an Euclidean space. In the pointed sense every linear transport
along paths  turns out to be locally (along any fixed path) Euclidean. (In
particular, this is true for for parallel transports generated by linear
connections.) Sect. 4 contains investigation of two equivalent descriptions of
the studied transports along paths. Firstly, the local description in terms of
their local coefficients, which are analogous to the ones of a linear
connection, and, secondly, the one in terms of derivations along paths in
vector bundles. At the end, in Sect. 5, some comments on the considered here
problems are made.

Many problems, such as the comparison with the general connection's theory, integrability conditions (connected with the concepts curvature and torsion) and application to mathematical or physical questions, remain out of the field of this paper and will be considered elsewhere.

\medskip
\medskip
 {\bf 2. DEFINITION AND GENERAL FORM OF LINEAR\\ TRANSPORTS
 				 ALONG PATHS IN VECTOR BUNDLES}

\medskip
 By $(E,\pi ,B)$ we denote a general real vector bundle [5,6], p.44. Here the
base $B$ and the total bundle space $E$ are topological spaces, the
projection $\pi :E  \to B$ maps $E$ onto $B$ and the structure of $E$ is such
that  the fibres $\pi ^{-1}(x)\subset E, x\in B$ are isomorphic real vector
spaces. By $J$ and $\gamma :J  \to B$ we denote, respectively, an arbitrary
real interval and a path in B.

Definition 2.1 of Ref. [3] describes the $S$-transports along paths in tensor bundles as maps having the mentioned in it properties. A simple overview on them reviles that part of them are specific for the considered in [3] tensor bundles, the other ones being independent of this fact. A deeper analysis of the above-cited definition shows that its straightforward generalization leads to the concept of a linear transport along paths in vector bundles which is fixed by

{\bf Definition} ${\bf 2}{\bf .}{\bf 1}{\bf .} A$ linear transport
$(L$-transport) along paths in the real vector bundle $(E,\pi ,B)$ is a map
$L$ which to any path $\gamma :J  \to B$ assigns a map $L^{\gamma },
L$-transport along $\gamma $, such that $L^{\gamma }:(s,t)\mapsto
L^{\gamma }_{s  \to t}$, where for every $s,t\in J$ the map
\[
 L^{\gamma }_{s  \to t}:\pi ^{-1}(\gamma (s))  \to \pi ^{-1}(\gamma (t)),
\qquad (2.1)
\]
called an $L$-transport along $\gamma $ from $s$ to $t$, has the following
three properties:
\[
L^{\gamma }_{s  \to t}(\lambda u+\mu v)=\lambda S^{\gamma }_{s  \to
t}u+\mu S^{\gamma }_{s  \to t}v,
\quad \lambda ,\mu \in {\Bbb R},\ u,v\in \pi^{-1}(\gamma (s)), \qquad (2.2)
\]
\[
 L^{\gamma }_{t  \to r}\circ L^{\gamma }_{s  \to t}
=L^{\gamma }_{s  \to r}, r,s,t\in J,\qquad (2.3)
\]
\[
L^{\gamma }_{s  \to s}={\it id}_{\pi ^{-1_{(\gamma (s))}}},\qquad (2.4)
\]
where {\it id}$_{X}$ means the identity map of the set X.

{\bf Remark.} If in this definition we admit $(E,\pi ,B)$ to be a complex
vector bundle [6] and $\lambda ,\mu \in {\Bbb C}$ (instead of $\lambda ,\mu
\in {\Bbb R})$, we get the definition of a (${\Bbb C}-)$linear transport
along paths in complex vector bundles. Almost all further definitions and
results are valid also in the complex case (and, generally, in a case of an
arbitrary field ${\Bbb K})$ but, for the sake of shortness, we shall
investigate only the real case.

Evident examples of $L$-transports along paths are the considered in [3]
restrictions of $S$-transports on concrete tensor bundles.

 Putting $r=s$ in (2.3) and using (2.4), we get
\[
\Bigl( L^{\gamma }_{s  \to t  } \Bigr)^{-1}
=L^{\gamma }_{t  \to s}, \quad s,t\in J.  \qquad (2.5)
\]
So, the linear transports along a path are linear isomorphisms of the fibres
over the points of that path.

The following two propositions establish the general form of the linear transports along paths.

{\bf Proposition} ${\bf 2}{\bf .}{\bf 1}{\bf .} A$ map (2.1) is a linear
transport along $\gamma $ from $s$ to $t$ for every $s,t\in J$ if and only if
there exist an isomorphic with $\pi ^{-1}(x), x\in B$ vector space $V$ and a
family of linear isomorphisms $\{F^{\gamma }_{s}:\pi ^{-1}(\gamma (s))  \to
V, s\in J\}$ such that
\[
 L^{\gamma }_{s  \to t}={\bigl(}F^{\gamma }_{t}{\bigr)}^{-1}\circ
{\bigl(}F^{\gamma }_{s}{\bigr)}, s,t\in J. \qquad (2.6)
\]

	{\bf Proof.} If (2.1) is an $L$-transport along $\gamma $ from $s$
to $t$, then fixing some $s_{0}\in J$ and using (2.3) and (2.5), we get
$L^{\gamma }_{s  \to t}=L^{\gamma }_{s_{0}}\circ L^{\gamma }_{s  \to
s_{0}}=\left(\begin{array}{c} \end{array}\right. L^{\gamma }_{t  \to
s_{0}}\circ L^{\gamma }_{s  \to s_{0}}$, i.e. (2.6) holds for $V=\pi
^{-1}(\gamma (s_{0}))$ and $F^{\gamma }_{s}=L^{\gamma }_{s  \to s_{0}}$. On
the contrary, if (2.6) is valid for some linear isomorphisms $F^{\gamma
}_{s}$, then a straightforward calculation shows that it converts (2.3) and
(2.4) into identities and (2.2) is true due to the linearity of $F^{\gamma
}_{s}$, i.e. $L^{\gamma }_{s  \to t}$is an $L$-transport along $\gamma $ from
$s$ to $t$ for every $s,t\in $J.\blacksquare

{\bf Proposition 2.2.} Let in the vector bundle $(E,\pi ,B)$ be defined an
$L$-transport along paths with a representation (2.6) for some vector space
$V$ and linear isomorphisms $F^{\gamma }_{s}:\pi ^{-1}(\gamma (s))  \to V,
s\in $J. Then for a vector space $\mathbf{ V}$ there exist linear isomorphisms
$\mathbf{ F}^{\gamma }_{s}:\pi ^{-1}(\gamma (s))  \to \mathbf{ V}, s\in J$ for which
\[
L^{\gamma }_{s  \to t}={\bigl(}\mathbf{ F}^{\gamma }_{t}
{\bigr)}^{-1}\circ {\bigl(}\mathbf{ F}^{\gamma }_{s}{\bigr)}, s,t\in J.
 \qquad (2.7)
\]
iff there exists a linear isomorphism $D^{\gamma }:V \to \mathbf{ V}$ such that
\[
\mathbf{ F}^{\gamma }_{s}=D^{\gamma }\circ F^{\gamma }_{s}, s\in J.
 \qquad (2.8)
\]

	 {\bf Proof.} If (2.8) holds, then substituting $F^{\gamma
}_{s}=(D^{\gamma })^{-1}\mathbf{ F}^{\gamma }_{s}$ into (2.6), we get (2.7) and
vice versa, if (2.7) is valid, then from its comparison with (2.6) follows
$D^{\gamma }:=(\mathbf{ F}^{\gamma }_{t})(F^{\gamma }_{t})^{-1}= =(\mathbf{
F}^{\gamma }_{s})(F^{\gamma }_{s})^{-1}$to be the needed (independent of $s$
and $t)$ isomorphism.\blacksquare

So, the linear transport along $\gamma $ form $s$ to $t$ decomposes into a composition of linear maps depending separately on $s$
and $t$ (cf. proposition 2.1), the arbitrariness of these maps being
described by $eq. (2.8)$ (see proposition 2.2).

For some purposes is useful the linear transports along paths to be described
by their matrix elements in corresponding bases (or frames), which we are going
to do now.

Let $\{e_{i}\}$, where here and henceforth the Latin indices run from 1 to
$\dim(\pi ^{-1}(x))$, $x\in B$ and the usual summation rule from 1 to
$\dim(\pi ^{-1}(x))$ over repeated on different levels indices is assumed, be a field
of bases along $\gamma :J  \to B$, i.e. for every $s\in J$ the set of vectors
$\{e_{i}(s)\}$ to be a basis in the vector space $\pi ^{-1}(\gamma (s))$.

Due to (2.1), we have $L^{\gamma }_{s  \to t}e_{i}(s)\in \pi ^{-1}(\gamma
(t))$, hence there exists a unique matrix $H(t,s;\gamma ):=
H^{i}_{.j}(t,s;\gamma )  $ such that
\[
L^{\gamma }_{s  \to t}e_{i}(s)=H^{j}_{.i}(t,s;\gamma )e_{j}(t), s,t\in J.
 \qquad (2.9)
\]
 From here immediately follows that $H^{i}_{.j}(t,s;\gamma )$  are  elements
of a (two-point) tensor from the tensor space $\pi ^{-1}(\gamma (t))\otimes
\otimes (\pi ^{-1}(\gamma (s)))^{*}$, where $\otimes $ is the tensor  product
sign  and  the asterisk means the dual of the corresponding vector space.

{\bf Definition 2.2.} The matrix function $H:(t,s;\gamma )\mapsto  \mapsto H(t,s;\gamma ), s,t\in J$ will be called a matrix of the considered $L$-transport along paths.

The matrix of the $L$-transport uniquely defines its action. In fact, if
$u=u^{i}e_{i}(s)\in \pi ^{-1}(\gamma (s))$, then from (2.2) and (2.9), we get
\[
 L^{\gamma }_{s  \to t}u=H^{i}_{.j}(t,s;\gamma )u^{j}e_{i}(t).\qquad (2.10)
\]
In terms of the matrix function $H$, due to (2.9), the basic properties (2.3)
and (2.4) take, respectively, the form
\[
H(r,t;\gamma )H(t,s;\gamma )=H(r,s;\gamma ), r,s,t\in J,\qquad (2.11)
\]
\[
 H(s,s;\gamma )={\bf 1} , s\in J,\qquad (2.12)
\]
 where ${\bf 1}$ is the unit matrix and in the matrix multiplication as a
first matrix index the superscript is understood.

On the contrary, due to (2.10), from (2.11) and (2.12) follow, respectively, (2.3) and (2.4). So, we have proved

{\bf Proposition} ${\bf 2}{\bf .}{\bf 3}{\bf .} A$ linear map (2.1) is a linear transport along $\gamma $ from $s$ to $t$ iff its matrix $H$, defined by (2.9), satisfies (2.11) and (2.12).

Let us introduce the matrices $F(s;\gamma ):=[F^{i}_{.j}(s;\gamma)]$ and
$\mathbf{ F}(s;\gamma ):=[\mathbf{ F}^{i}_{.j}(s;\gamma )] $, where
$F^{i}_{.j}(s;\gamma )$ and $\mathbf{ F}^{i}_{.j}(s;\gamma )$ are the matrix
elements, respectively, of the maps $F^{\gamma }_{s}$ and $\mathbf{ F}^{\gamma
}_{s}$, i.e. $F^{\gamma }_{s}e_{j}(s)=:F^{i}_{.j}(s;\gamma )l_{i}$ and $\mathbf{
F}^{\gamma }_{s}e_{j}(s)=:\mathbf{ F}^{i}_{.j}(s;\gamma )\mathbf{ l}_{i}$ in which
$\{l_{i}\}$ and $\{\mathbf{ l}_{i}\}$ are, respectively, bases in $V$ and $\mathbf{
V}$. Then, using (2.2) and (2.9), we see that in terms of matrix elements
propositions 2.1 and 2.2 are, respectively, equivalent to

{\bf Proposition 2.4.} The linear map (2.1) with defined by (2.9) matrix
elements in some frames is an $L$-transport along $\gamma $ from $s$ to $t$
iff there exist nondegenerate matrices $F(s;\gamma ), s\in J$ such that
\[
H(t,s;\gamma )=(F(t;\gamma ))^{-1}F(s;\gamma ), s,t\in J.
(2.13)
\]
 {\bf Proposition 2.5.} Let the matrix $H$ of some $L$-transport along
paths has the representation (2.13). Then for certain nondegenerate matrices
$\mathbf{ F}(s;\gamma ), s\in J$ is valid the equality
\[
H(t,s;\gamma )=(\mathbf{ F}(t;\gamma ))^{-1}\mathbf{ F}(s;\gamma ), s,t\in
J\qquad (2.13^\prime )
\]
 if and only if there exists a nondegenerate matrix $D(\gamma )$, depending
only on $\gamma $, such that
\[
\mathbf{ F}(s;\gamma )=D(\gamma )F(s;\gamma ), \quad s\in J,\
\det D(\gamma )\neq 0,\infty .\qquad (2.14)
\]

\medskip
 {\bf 3. SPECIAL BASES FOR LINEAR\\ TRANSPORTS ALONG PATHS}

\medskip
Special frames for an $L$-transport along paths we call frames in which it
(generally locally; see below) looks like a parallel transport in an Euclidean
spaces. The existence of such frames is expressed by

{\bf Proposition 3.1.} For every $L$-transport along a path $\gamma $ there
exist a basis in which the transport's matrix is unit{\bf .} If in a given
basis along a path $\gamma $ the matrix of some $L$-transport is constant or
depends only on $\gamma $, then it is unit{\bf .} All frames along $\gamma $ in
which the matrix of an $L$-transport is unit are obtained from the above one,
and, consequently, from one another, by a nondegenerate linear transformations
with constant or depending on $\gamma $ coefficients.

As from the view point of the applications this proposition is an important one, we shall give two its proofs. The first proof admits a straightforward application for different purposes, and the second one distinguishes by its shortness and clearness.

{\bf Proof 1.} Let $\{e_{i}(s)\}$ be an arbitrary basis in $\pi ^{-1}(\gamma
(s))$ and $\{e^{i}(s)\}$ be its dual basis, $s\in $J. By proposition 2.4
there exist nondegenerate matrices $F(s;\gamma ), s\in J$, such that
$H^{i}_{.j}(t,s;\gamma )=((F(s;\gamma ))^{-1})^{i}_{.k}(F(s;\gamma
))^{k}_{.j}  $ If the frames $\{e^{i^\prime }(s)\}, s\in J$ are defined
through
\[
 e^{i^\prime }(s):=\delta ^{i^\prime }_{l}(F(s;\gamma
))^{l}_{.i}e^{i}(s)),\qquad (3.1)
\]
 where $\delta ^{i}_{l}$are the Kronecker  deltas $(\delta ^{i}_{l}=1$ for
$i=j$ and $\delta ^{i}_{l}=0$ for $i\neq j)$, then
 \[
 e_{j^\prime }=\delta ^{m}_{j^\prime }((F(s;\gamma
))^{-1})^{j}_{.m}e_{j}\qquad (3.2)
\]
 and due to the fact that $H^{i}_{.j}(t,s;\gamma )$ are components of a
tensor from $\pi ^{-1}(\gamma (t))\otimes (\pi ^{-1}(\gamma (s)))^{*}$, we
get $H^{i^\prime }_{.}(t,s;\gamma )= =[\delta ^{i^\prime }_{l}(F(t;\gamma
))^{l}_{.i}][\delta ^{m}_{j^\prime }((F(s;\gamma
))^{-1})^{j}_{.m}]H^{i}_{.j}(t,s;\gamma )=\delta ^{i}_{j}$, i.e. in the basis
$\{e_{i^\prime }\}$ the matrix of the considered $L$-transport along $\gamma
$ is unit.

Let now in $\{e_{i}\}$ to be fulfilled $H(t,s;\gamma )=A(\gamma ), A(\gamma
)$ being a nondegenerate matrix depending possibly on the path $\gamma $.
According to (2.13), we have
$A(\gamma )=(F(t;\gamma ))^{-1}F(s;\gamma )$ and
consequently $F(s;\gamma )=A(\gamma )F(t;\gamma )=:B(\gamma )$, where, due to
the arbitrariness of $s$ and $t, B(\gamma )$ depends only on $\gamma $. Hence
$H(t,s;\gamma )=(B(\gamma ))^{-1}B(\gamma )={\bf 1}$ and taking into account
the above definition of $\{e_{i^\prime }(s)\}$, we get $e_{j}=(B(\gamma
))^{i}_{.j}\delta ^{i^\prime }_{i}e_{i^\prime }$which shows that all frames in
which the matrix of the $L$-transport is constant (unit) are obtained from
$\{e_{i^\prime }(s)\}$ by linear transformations with constant or depending
on the path coefficients, so they themselves are connected with such
transformations. To end this proof it remains only to be noted that due to
(2.14) the dependence of $B(\gamma )$ on $\gamma $ is insignificant and it
may be considered only the case $B(\gamma )=$const, which is equivalent to
the not changing $H(s,t;\gamma )$ redefining $F(s;\gamma ).\blacksquare $

 {\bf Proof 2.} Let us fix $s_{0}\in J$ and a basis $\{f'\}$ in
$\pi ^{-1}(\gamma (s_{0}))$. Along $\gamma :J  \to B$ we define a basis
$\{e'_i\}$ such that in $\pi ^{-1}(\gamma (s))$ it has the form
\[
e'_i(s):=I^{\gamma }_{s_{0}\to s} f'_i  \qquad (3.3)
\]
Then
\(
I^{\gamma }_{s  \to t}e'_i(s)
=I^{\gamma }_{s  \to t}\circ I^{\gamma}_{s_{0}\to s}f'_i
=I^{\gamma }_{s_{0}\to t}f'_i=e'_i(t),
\)
 hence, due to (2.9), the matrix of $L$ in $\{e'\}$ is
$H^\prime(t,s;\gamma )={\bf 1}$.

Let $\{e^{\gamma }_{i}\}$ be arbitrary basis along $\gamma $ in which the matrix
of $L$ is $H(t,s;\gamma )=B(\gamma )$ for some, may be depending on $\gamma$,
nondegenerate matrix $B(\gamma )$. Then there exists a nondegenerate matrix
$ [A^{j}_{i}(s;\gamma ) ]$ such that
$e^{\gamma }_{i}(s)=A^{j}_{i}(s;\gamma )e'_j$ (s). The substitution of the
last equality into $I^{\gamma }_{s  \to t}e^{\gamma }_{i}(s)=(B(\gamma
))^{j}_{i}e^{\gamma }_{j}(t) ($see (2.9)), due to (2.10), gives
$A^{j}_{i}(s;\gamma )=(B(\gamma ))^{k}_{i}A^{j}_{k}(t;\gamma )$. From here
for $t=s$, we get $B(\gamma )={\bf 1}$, so the same equality reduces to
$A^{j}_{i}(s;\gamma )=A^{j}_{i}(t;\gamma )$ and consequently
$A^{j}_{i}(s;\gamma )$ are constants or may depend only on $\gamma
.\blacksquare $

In connection with the fact that generally the global version of proposition
3.1, which will be considered elsewhere, is not true it is convenient to be
introduced the concept of an Euclidean case concerning $L$-transports along
paths, i.e. when in a vector bundle $(E,\pi ,B)$ is given an Euclidean
$(L-)$transport along paths defined by

{\bf Definition 3.1.} An $L$-transport
in a vector bundle $(E,\pi ,B)$ is Euclidean over the set $U\subset B$, or it
is Euclidean if $U=B$, if in $\pi ^{-1}(U)$ there exists a field of bases
$\{e_{i}\}$, i.e. $\{e_{i}(x)\}$ is a basis in $\pi ^{-1}(x)$ for $x\in U$,
in which the matrix of the transport is unit along any path in $U$, i.e. if
$\gamma :J  \to B$ and there is an ${\Bbb R}$-interval $J^\prime \subset J$
such that $\gamma (J^\prime )\subset U$, then $H(t,s;\gamma )={\bf 1}$, for
$s,t\in J^\prime $.

{\bf Definition 3.2.} Let $\{e_{i}\}$ be a basis in a
vector bundle $(E,\pi ,B)$, i.e. $\{e_{i}(x)\}$ is a basis in $\pi ^{-1}(x),
x\in $B. A generated by (associated to$) \{e_{i}\}$ Euclidean $(L-)$transport
in $(E,\pi ,B)$ is an $L$-transport in this fibre bundle the matrix of which
in the basis $\{e_{i}\}$ along any path is unit, i.e. $H^{i}_{.j}(s,t;\gamma
)=\delta ^{i}_{j}$ for every path $\gamma :J\to B$ and every $s,t\in $J.

The name "Euclidean transport" is connected with the fact that if above we put $B={\Bbb R}^{n}$and identify $T_{x}({\Bbb R}^{n})$ and ${\Bbb R}^{n}$, then in any orthogonal basis, i.e. in Cartesian coordinates,
the Euclidean transport coincides with the standard parallel transport in ${\Bbb R}^{n}$under which the vector's components are left unchanged.

{\bf Proposition 3.2.} Every basis in a vector bundle generates a unique Euclidean transport in it. The opposite statement being valid only locally, i.$e$ along a given path.

{\bf Proof.} The first path of the proposition is a corollary of proposition 2.4 and definition 2.2.

By proposition 3.1 for every $L$-transport along a fixed path there exists a basis along it in the bundle in which the transport's matrix is (locally) unit and, hence, the generated from this basis Euclidean transport along the pointed path coincides with the initial $L$-transport.\blacksquare

The importance of the last proposition is in the fact that in the above
sense \emph{any ${L}$-transport is locally Euclidean}.

{\bf Proposition 3.3.} Two (or more) frames generate one and the same Euclidean
transport in a given vector bundle if and only if they can be obtained from
each other through linear transformations with constant coefficients.

{\bf Proof.} This result is direct consequence from proposition 2.4 and definition 3.1.\blacksquare

The above results show that any $L$-transport is Euclidean over any point of
the base (see (2.12)), as well as over an arbitrary path in it (see
proposition 3.1). For other subsets $U\subset B$ this is, generally, not
true. In particular, for $U=B$ analogous result holds only for "flat$"
L$-transports, whose curvature operator vanishes, a result which will be
establish in another paper.

\medskip
\medskip
 {\bf 4. THE EQUIVALENCE OF LINEAR TRANSPORTS\\ ALONG PATHS
	AND DERIVATIONS ALONG PATHS}

\medskip
Let a linear transport along paths $L$ in the vector bundle $\xi =(E,\pi ,B)$ be given. Let it to be smooth of class $C^{1}$in a sense that such is its matrix $H$ with respect to its first, and consequently to its second (see (2.13)), argument. Let $\xi \mid _{\gamma (J)}:=(\pi ^{-1}(\gamma (J)),\pi \mid _{\gamma (J)},\gamma (J))$ be the restriction of $\xi $ on the set $\gamma (J)$ defined by the path $\gamma :J  \to $B. Let Sec$(\xi ) ($resp. Sec$^{k}(\xi ))$ be the set of all (resp. of class $C^{k})$ sections of $\xi  [6,1]$.

 {\bf Definition 4.1.} The derivation along paths generated by $L$ is a map
${\cal D}$ such that ${\cal D}:\gamma \mapsto {\cal D}^{\gamma }$,
where the assigned to any path $\gamma $ map
\[
 {\cal D}^{\gamma }:Sec^{1}(\xi \mid _{\gamma (J)})  \to Sec(\xi \mid
_{\gamma (J)}),\qquad (4.1)
\]
 called derivation along $\gamma $ generated by
$L$, for every $\sigma \in $Sec$^{1}(\xi \mid _{\gamma (J)})$ is defined by
${\cal D}^{\gamma }:\sigma \mapsto {\cal D}^{\gamma }\sigma $ in
which for $s,s+\epsilon \in J$
\[
{\bigl(}{\cal D}^{\gamma }\sigma {\bigr)}(\gamma (s))
:={\cal D}^{\gamma }_{s}\sigma
:=\lim_{\varepsilon\to 0}\bigl[ \frac{1}{\varepsilon}
{\bigl(} L^{\gamma }_{s+\epsilon \to s}\sigma (\gamma (s+\epsilon ))
-\sigma (\gamma (s)){\bigr)}\bigr].   \qquad (4.2)
\]

	The derivative of $\sigma $ along $\gamma $ with respect to $L$ is
${\cal D}^{\gamma }\sigma $ whose value at $\gamma (s)$, i.e. $({\cal
D}^{\gamma }\sigma )(\gamma (s))$, is determined by the map
\[
{\cal D}^{\gamma }_{s}:Sec^{1}(\xi \mid _{\gamma (J)})  \to \pi
^{-1}(\gamma (s)).\qquad (4.3)
\]
 The limit in (4.2), due to the above
condition of smoothness, always exist.

 {\bf Proposition 4.1.} The derivation ${\cal D}^{\gamma }$ is an
${\Bbb R}$-linear map, i.e. for
$\lambda _{1},\lambda _{2}\in {\Bbb R}$ and
$\sigma_{1},\sigma _{2}\in Sec^{1}(\xi \mid _{\gamma (J)})$ is fulfilled
\[
 {\cal D}^{\gamma }(\lambda _{1}\sigma _{1}+\lambda _{2}\sigma
_{2})=\lambda _{1}{\cal D}^{\gamma }\sigma _{1}+\lambda _{2}{\cal D}^{\gamma
}\sigma _{2},\qquad (4.4)
\]
 and its value at every $t\in J$ satisfies the identity
\[
{\cal D}^{\gamma }_{t}\circ L^{\gamma }_{s  \to t}\equiv 0, s,t\in J.
 \qquad (4.5)
\]

{\bf Proof.} (4.4) and (4.5) follow from (4.2) and, respectively, (2.2)
and (2.3).\blacksquare

 Let along $\gamma :J  \to B$ be given a field of bases $\{e_{i}\}$, i.e.
$\{e_{i}(s)\}$ to be a basis in $\pi ^{-1}(\gamma (s))$, in which the
$L$-transport along paths $L$ to be defined by its matrix $H:(t,s;\gamma
)\mapsto  \mapsto H(t,s;\gamma ):=  H^{i}_{.j}(t,s;\gamma )  $. Then along
$\gamma $ every $\sigma \in $Sec$(\xi \mid _{\gamma (J)})$ has a unique
representation in the form $\sigma =\sigma ^{i}e_{i}$, where $\sigma
^{i}:\gamma (J)  \to {\Bbb R}$ are the components of $\sigma $ in
$\{e_{i}\}$. Hence $\sigma (\gamma (s))= =\sigma ^{i}(\gamma (s))e_{i}(s)$
and from (4.2) and (2.10), we get
\[
 {\bigl(}{\cal D}^{\gamma }\sigma {\bigr)}(\gamma (s)):={\cal
D}^{\gamma }_{s}\sigma ={\bigl(}{\cal D}^{\gamma }_{s}\sigma
{\bigr)}^{i}e_{i}(s)=
\]
\[
=\bigl\{ \frac{\partial}{\partial\varepsilon}
\bigl[ \sum_j H^{\hbox{i.}}_{.j}(s,s+\epsilon
;\gamma )\sigma ^{j}(\gamma (s+\epsilon ))\bigr]{\bigr\}}
|_{\epsilon =0}e_{i}(s).\qquad (4.6)
\]
 From here immediately follows

{\bf Proposition 4.2.} The explicit action of ${\cal D}^{\gamma }$on $\sigma
\in  \in $Sec$^{1}\xi \mid _{\gamma (J)})$ is
\[
{\bigl(}{\cal D}^{\gamma }\sigma {\bigr)}(\gamma (s))
=
\Bigl[
\frac{d \sigma^i(\gamma(s))}{ds}
+ \Gamma ^{i}_{.j}(s;\gamma )\sigma^{j}(\gamma (s))
\Bigr] e_{i}(s),  \qquad (4.7)
\]
where
\[
 \Gamma ^{i}_{.j}(s;\gamma )
:= \frac{\partial H_{.j}^{i}(s,t;\gamma)}{\partial t}\Big|_{t=s}
= - \frac{\partial H_{.j}^{i}(t,s;\gamma)}{\partial t}\Big|_{t=s}
\qquad (4.8)
\]
(The last equality follows, e.g., from (2.13).)

It is convenient in any basis to introduce the matrix function
$\Gamma_{\gamma }:s\mapsto \Gamma _{\gamma }(s)
:=[\Gamma ^{i}_{.j}(s;\gamma )],\ s\in $J.
If the matrix $H$ of the linear transport along paths has a
representation (2.13) (see proposition 2.4), then, due to (4.8), we obtain
\[
\Gamma _{\gamma }(s)
= \frac{\partial H(s,t;\gamma)}{\partial t}\Big|_{t=s}
=F^{-1}(s;\gamma ) \frac{d F(s;\gamma)}{d s}
 \qquad (4.9)
\]
As we shell prove (see below proposition 4.5) the quantities (4.8) give
an adequate description of the transport L. They also have a sense of a
"depending on the path $\gamma $ coefficients of a linear connection$" (cf.
[6,1])$ which is confirmed by

{\bf Proposition 4.3.} Let the basis $\{e_{i}(s)\}$ in $\pi ^{-1}(\gamma
(s))$ be changed to $\{e_{i^\prime }(s)=A^{i}_{i^\prime }(s)e_{i}(s)\},
A(s):= [A^{i}_{i^\prime }(s)]$ being a nondegenerate matrix and
$A^{-1}(s):=[A^{i^\prime }_{i}(s)]$.
The transformation $\{e_{i}\}  \to \{e_{i^\prime
}=A^{i}_{i^\prime }e_{i}\}$ leads to the transformation of $H$ and $\Gamma
_{\gamma }$, respectively, into
\[
 H^\prime (s,t;\gamma )=A^{-1}(s)H(s,t;\gamma )A(t),\qquad (4.10)
\]
\[
\Gamma ^\prime (s;\gamma )
=A^{-1}(s)\Gamma (s;\gamma )A(s)+A(s)^{-1}  \frac{d A(s)}{ds} ,\qquad (4.11)
\]
which in component form, respectively, read
\[
H^{i^\prime }_{..j^\prime }(s,t;\gamma )=A^{i^\prime
}_{i}(s)A^{j}_{j^\prime }(t)H^{i}_{.j}(s,t;\gamma ),\qquad (4.10^\prime )
\]
\[
\Gamma ^{i^\prime }_{..j^\prime }(s;\gamma )
=A^{i^\prime}_{i}(s)A^{j}_{j^\prime }(s)\Gamma ^{i}_{.j}(s;\gamma )
+A^{i^\prime }_{i}(s)  \frac{d A_{j^\prime }^{i}(s)}{ds}
 \qquad (4.11^\prime )
\]

 {\bf Proof.} $Eq.  (4.10)$ expresses the fact that
$H^{i}_{.j}(s,t;\gamma )$ are components of a tensor from $\pi ^{-1}(\gamma
(s))\otimes (\pi ^{-1}(\gamma (t)))^{*}$and $eq. (4.11)$ is a corollary of
$(4.10), (2.12)$ and (4.9). The equivalence of (4.10) and $(4.10^\prime )$
and (4.11) and $(4.11^\prime )$ follows from the definitions of the
corresponding matrices.\blacksquare

So, any (smooth$) L$-transport along paths generates along every path $\gamma $ functions (4.8) transforming according to $(4.11^\prime )$. Here, naturally, arises the opposite problem. Let along every path $\gamma $ in any basis $\{e_{i}\}$ along it be defined functions $\Gamma ^{i}_{.j}(s;\gamma )$ which, when the basis is changed, transform
through $(4.11^\prime )$. Does there exist an $L$-transport along paths generating these functions by means of $eq. (4.8)$? The answer of this question is given by

 {\bf Proposition 4.4.} Let in any basis $\{e_{i}\}$ along a path $\gamma :J
\to B$ be given (continues with respect to $s)$ functions $\Gamma
^{i}_{.j}(s;\gamma )$ which, when $\{e_{i}\}$ is changed, have the
transformation low $(4.11^\prime )$ and $\Gamma _{\gamma }(s):=[\Gamma
^{i}_{.j}(s;\gamma )]$. Then there exists a unique $L$-transport along paths
generating along $\gamma $ the matrix $\Gamma _{\gamma }$ by (4.9) and the
matrix of which, $H:($s.$t;\gamma )\mapsto H($s.$t;\gamma )= =F^{-1}(s;\gamma
)F(t;\gamma )$, is
\[
 H(s,t;\gamma )
=Y(s,s_{0};-\Gamma _{\gamma })[Y(t,s_{0};-\Gamma _{\gamma })]^{-1},
\quad s,t\in J,\qquad (4.12)
\]
 where $s_{0}\in J$ is fixed and for any
continues matrix function $Z:s\mapsto Z(s)$ the matrix $Y:=  Y_{ij}
:=Y($s.$s_{0};Z)$ is the unique solution of the initial-value problem
\[
\frac{dY}{ds} = Z(s)Y,\quad Y=Y(s,s_{0};Z),\quad s\in J,\qquad (4.13a)
\]
\[
Y(s_{0},s_{0};Z)= {\bf 1} .\qquad (4.13b)
\]

 {\bf Proof.} Let us note at the beginning that the existence and uniqueness
of the solution of (4.13) is proved in [7].

At first we shall prove that (4.12) is the unique solution of (4.9) when $\Gamma _{\gamma }$is given. In fact, using $H($s.$t;\gamma )= =F^{-1}(s;\gamma )F(t;\gamma ) ($see proposition 2.4) and $dF^{-1}/ds= =-F^{-1}(dF/ds)F^{-1}[8]$, we see $eq. (4.9)$ with respect to $H$ to be equivalent to the equation $dF^{-1}(s;\gamma )/ds=-\Gamma _{\gamma }(s;\gamma )F^{-1}(s;\gamma )$. Arbitrary fixing $s_{0}\in J$ and comparing this equation (with respect to $F^{-1})$ with (4.13), we see that its solutions is $F^{-1}(s;\gamma )=Y(s,s_{0};-\Gamma _{\gamma })$. Substituting this expression for $F^{-1}$into (2.13), we get (4.12), which is independent of the concrete choice of $F(s_{0};\gamma )$ and $s_{0}($see proposition 2.5 and the
given in [7] properties of $Y)$.

So, in any fixed basis only the matrix (4.12) generates the given $\Gamma _{\gamma }$by (4.9). Hence, defining an $L$-transport along paths which is such a basis has a matrix (4.12), we see, in conformity with (4.10) and (4.11), that in every basis along $\gamma $ it generates the functions $\Gamma ^{i}_{.j}(s;\gamma )$ through (4.8). By construction this $L$-transport along paths is unique.\blacksquare

 From propositions 4.3 and 4.4 directly follows

{\bf Proposition 4.5.} Let in any basis along every path $\gamma :J  \to B$ be given a set of functions $\{\Gamma ^{i}_{.j}(s;\gamma )\}$. Then, when the basis is changed, they transform in accordance with $(4.11^\prime )$ iff there exists a (unique) linear transport along paths generating them through its matrix by (4.8).

Consequently, the definition of an $L$-transport along a path $\gamma $ is equivalent to the definition in any basis along $\gamma $ of functions $\Gamma ^{i}_{.j}(s;\gamma )$ having the transformation low $(4.11^\prime )$. This is a reason to give

{\bf Definition 4.2.} The functions $\Gamma ^{i}_{.j}:(s;\gamma )\mapsto \Gamma ^{i}_{.j}(s;\gamma )$ assigned to any linear transport along paths through (4.8) or defining such a transport by means of (4.12) will be called coefficients of that $L$-transport along paths.

In proposition 4.2 we saw that the coefficients $\Gamma ^{i}_{.j}$of an
$L$-transport along paths uniquely define the action of the operator (4.1).
Besides, if $e_{i}\in $Sec$(\xi \mid _{\gamma (J)})$ form a basis along
$\gamma :J  \to B$, i.e. if $\{e_{i}(s)\}$ is a basis in $\pi ^{-1}(\gamma
(s))$, then (4.7) results in
\[
 {\cal D}^{\gamma }e_{j}=(\Gamma _{\gamma })^{i}_{.j}e_{i},\qquad (4.14)
\]
which is equivalent to
\[
{\cal D}^{\gamma }_{s}e_{j}
=\Gamma ^{i}_{.j}(s;\gamma )e_{i}(s)\qquad (4.14^\prime )
\]
and can be used as an equivalent to (4.8) definition of the coefficients
$\Gamma ^{i}_{.j}(s;\gamma )$.

As the coefficients $\Gamma ^{i}_{.j}(s;\gamma )$ give an adequate description of the $L$-transports along paths, the $eq. (4.14)$ suggests that such a description can be done also in terms of operators (4.1) with an "appropriate properties". Before formulating the corresponding result we need to give one definition.

If $\sigma \in $Sec$^{1}(\xi \mid _{\gamma (J)})$ and $f:J  \to {\Bbb R}$ is
$a C^{1}$ function, then from (4.7), we find
\[
{\cal D}^{\gamma }_{s}(f\cdot \sigma )
=  \frac{df(s)}{ds} \cdot \sigma (\gamma (s))+f(s)\cdot ({\cal D}^{\gamma
}_{s}\sigma ).\qquad (4.15)
\]
If $B$ is a differentiable manifold, $\gamma $  is $a C^{1}$path and is the
tangent to $\gamma $ vector field, then (4.15) is equivalent to
\[
 {\cal D}^{\gamma }(f\cdot \sigma )=(\dot\gamma(f))\cdot \sigma +f\cdot ({\cal
D}^{\gamma }\sigma ).\qquad (4.16)
\]

 {\bf Definition} ${\bf 4}{\bf .}{\bf 3}{\bf .} A$ derivation along a path
$\gamma  ($in the vector bundle $\xi =(E,\pi ,B))$ is a linear map (6.1)
satisfying (4.15) in which ${\cal D}^{\gamma }_{s}$is defined by ${\cal
D}^{\gamma }_{s}\sigma :=({\cal D}^{\gamma }\sigma )(\gamma (s))$ for $\sigma
\in  \in $Sec$^{1}(\xi \mid _{\gamma (J)}). A$ derivation along paths is a
map assigning to any path a derivation along it.

Evident example of derivations along a path are the derivations of the tensor algebra (over a manifold) restricted over that path [3,1], as well as the maps (4.1) generated by $L$-transports along paths by means of (4.2).

{\bf Proposition 4.6.} The map ${\cal D}^{\gamma }$:Sec$^{1}(\xi \mid _{\gamma (J)})  \to $Sec$(\xi \mid _{\gamma (J)})$ is a derivation (in a sense of definition 4.3) iff there exists an $L$-transport $L^{\gamma }$along $\gamma $ generating it through (6.2).

{\bf Proof.} If ${\cal D}^{\gamma }$is generated through (6.2) by some $L$-transport along $\gamma $, then, as we proved, it satisfies (4.4) and (4.5), so it is a derivation along $\gamma $. Vice versa, let ${\cal D}^{\gamma }$be a
derivation along $\gamma $. Defining the functions $\Gamma ^{i}_{.j}(s;\gamma )$ through $(4.14^\prime )$ in any basis along $\gamma $, we find, due to (4.4) and (4.15), the $eq. (4.7)$ to be valid in any basis $\{e_{i}\}$ along $\gamma $. Thus we see that in the basis $\{e_{i^\prime }=A^{i}_{i^\prime }e_{i}\}$ is fulfilled $\Gamma ^{i^\prime }_{..j^\prime }(s;\gamma )e_{i^\prime }(s)={\cal D}^{\gamma }_{s}(e_{j^\prime })={\cal D}^{\gamma }_{s}(A^{j}_{j^\prime }e_{j})$ and using (4.4) and (4.15), we get that $\Gamma ^{i^\prime }_{..j^\prime }$and $\Gamma ^{i}_{.j}$are connected by $eq. (4.11^\prime )$. Hence, $\Gamma ^{i}_{.j}$are coefficients of some $L$-transport along paths $L$, which by proposition 4.5  exists, is unique, and generates them through (4.8). Besides, by proposition 4.2, the action of the assigned to it operator (4.1) is given also by (4.7), i.e. that operator, given explicitly by (4.2), coincides with ${\cal D}^{\gamma }.\blacksquare $

In the above proof of proposition 4.6 we saw that for any derivation ${\cal D}^{\gamma }$along $\gamma $ its action is uniquely defined by the functions $\Gamma ^{i}_{.j}$, determined uniquely by $(4.14^\prime )$. That is why we call them {\it coefficients} of ${\cal D}^{\gamma }$. From proposition 4.4 and the same proof immediately follows

{\bf Proposition 4.7.} Every $L$-transport $L^{\gamma }$along $\gamma $ defines through $(4.2) a$ unique derivation along $\gamma $ whose explicit action is given by (4.7) and whose coefficients coincide with those of $L^{\gamma }$and, vice versa, any derivation ${\cal D}^{\gamma }$along $\gamma $ defines a unique $L$-transports along $\gamma $ who generates ${\cal D}^{\gamma }$by (4.2), whose coefficients coincide with those of ${\cal D}^{\gamma }$and whose matrix is given by (4.12).

Practically, the last proposition establishes the equivalence between the set of all linear transports along paths and the one of all derivations along paths.

\medskip
\medskip
 {\bf 5. CONCLUDING REMARKS}

\medskip
We defined the linear transports along paths in vector bundles as maps satisfying part of the basic axioms of the $S$-transports (parallel transports generated by derivations of tensor algebras) [3]. The ties of the $L$-transports along paths with the connection and parallel transports theories will be examined from a more general view-point elsewhere. Here we want only to present these connections in the special case of linear connection and the corresponding to them parallel transports.

Let the base of the vector bundle $(E,\pi ,B)$, i.e. $B$, be a
differentiable manifold and $\gamma :J  \to B$ be $a C^{1}$path in it with a
tangent vector field   . Let us consider those $L$-transports along paths
whose coefficients in some basis along $\gamma $ have the form
\[
 \Gamma ^{i}_{.j}(s;\gamma )
= \sum_{\alpha=1}^{\dim (B)} \Gamma ^{i}_{.j\alpha }(\gamma (s))
\dot\gamma^{\alpha }(s),\qquad (5.1)
\]
 where $\Gamma ^{i}_{.j\alpha }:\gamma
(J)\to {\Bbb R}, i,j=1,\ldots  ,\dim(\pi ^{-1}(x)), x\in B, \alpha =
=1,\ldots  ,\dim(B)$ are some real functions and   $\dot\gamma^{\alpha }(s)$
are the local components of   (s) in some basis in the tangent to $B$ space
at $\gamma (s)$.

{\bf Proposition 5.1.} If the representation (5.1) holds in one basis
$\{e_{i}\}$ along $\gamma $, then it is true in a basis $\{e_{i^\prime
}=A^{i}_{i^\prime }e_{i}\}$ along $\gamma $ iff under the change of the frames
$\{e_{i}\}  \to \{e_{i^\prime }\}$ and of the coordinates $\{x^{\alpha }\}
\to \{x^{\alpha ^\prime }\}$ in a neighborhood of $x=\gamma (s)\in B$ the
functions $\Gamma ^{i}_{.j\alpha }(x)$ transform into
\[
\Gamma ^{i^\prime }_{..j^\prime \alpha ^\prime }(x)
= \sum_{i,j=1}^{\dim(\pi^{-1}(x))}
\Bigl( \sum_{\alpha=1}^{\dim(B)}
	A^{i^\prime }_{i}(s)A^{j}_{j^\prime }(s)
\frac{\partial x^\alpha}{\partial x^{\alpha'}}\Big|_x
\Gamma ^{i}_{.j\alpha }(x) \Bigr) +
\]
\[
+ \sum_{i=1}^{\dim(\pi^{-1} (x)) } A^{i^\prime }_{i}(s)
\frac{\partial A_{j^\prime }^{i}(s)}{\partial x^{\alpha'}}\Big|_x ,\quad
x=\gamma (s)\in B\qquad (5.2)
\]
 the right-hand-side of which depends on the map $\gamma $ only through the
point $x=\gamma (s)$.

{\bf Proof.} The proposition is a consequence of the substitution of (5.1) (in $\{e_{i}\}$ and in $\{e_{i^\prime }\})$ into $(4.11^\prime ).\blacksquare $

The transformation low (5.2) is a straightforward generalization of the one of the coefficients of a linear connection [1,2].

In particular, if $(E,\pi ,B)$ is the tangent bundle to the manifold $B$, i.e. $E=T(B)$ and $\pi ^{-1}(x)=T_{x}(B), x\in B$, then $\dim(B)= =\dim(\pi ^{-1}(x)), x\in B$ and putting $A^{i^\prime }_{i}(s)=\partial x^{i^\prime }/\partial x^{i}\mid _{x=\gamma (s)}$in (5.2), we see that $\Gamma ^{i}_{.j\alpha }(x), i,j,\alpha =1,\ldots  ,\dim(B)$ are coefficients of a linear connection (covariant differentiation$) \nabla  [1]$. In this case, due to (4.7), the assigned to the transport operator ${\cal D}^{\gamma }$is simply $\nabla _{\cdot }$, the covariant differentiation along $\gamma $, and hence the transport itself coincides with the corresponding to $\nabla $ parallel transport $[1] (cf. [4]$, Sect.5).

 In Sect. 3 we proved that for any liner transport along a path $\gamma $ in
a vector bundle there is along $\gamma  a$ class of special frames in which
the transport's matrix is unit, i.$e H(s,t;\gamma )={\bf 1}$. In every such
basis the transport's coefficients vanish as in them (see (4.9))
\[
 \Gamma _{\gamma }(s)=\partial {\bf 1}/\partial s=0\qquad (5.3)
\]
 and vice versa, any basis in which (5.3) holds belongs to this class.

 So, equation (4.7) in these frames is equivalent to
\[
 ({\cal D}^{\gamma }\sigma )(\gamma (s))
=   \frac{d\sigma^i(\gamma(s))}{ds} e_{i}(s),\qquad (5.4)
\]
i.e. ${\cal D}^{\gamma }$ behaves like a derivative along $\gamma $ of
vectors in an Euclidean space (in Cartesian coordinates).

 Thus we proved

{\bf Proposition 5.2.} In a given basis along a fixed path the following
three statements are equivalent:

 a) The matrix of the $L$-transport is unit.

 b) The coefficients of the $L$-transport vanish.

 c) The action of the assigned to the $L$-transport derivation reduces to a
simple derivation of the components of the sections with respect to the
path's parameter, i.e. (5.4) holds.

If in a vector bundle over a manifold the coefficients of an $L$-transport
along paths have the form (5.1) in some of the above special along $\gamma $
frames, then (5.3) leads to
\[
\sum_{\alpha=1}^{\dim(B)}\Gamma ^{i}_{.j\alpha }(\gamma (s))
\dot\gamma^{\alpha }(s)=0,\quad
 i,j=1,\ldots  ,\dim(\pi ^{-1}(\gamma (s))). (5.5)
\]
 But, nevertheless of the arbitrariness of $\gamma $, from here does
not follow $\Gamma ^{i}_{.j\alpha }(\gamma (s))=0$ as the basis itself,
generally, depends on the path $\gamma $.

In this connection let us note that frames in which (5.3), or in some more
special cases (5.5), is valid are a far going generalization of the coordinates
(or frames) in which the components of a symmetric affine connection vanish [2].

Let us also mention that to operations with vector bundles endowed with transports along paths there correspond analogous operations with these transports. In particular, to the direct sum and the tensor product of vector bundles [6] correspond, respectively, the direct sum and the tensor product of the linear transports in them, the matrices of which
are, respectively, the direct sum and the tensor product of the matrices of the initial transports along paths.

More concrete properties of the linear transports along paths in vector bundles, such a their curvature an torsion, as well as applications will be a subject of other works.

\medskip
\medskip
 {\bf ACKNOWLEDGEMENTS}

\medskip
The author expresses his gratitude to Prof. Vl. Aleksandrov (Institute of Mathematics of Bulgarian  Academy of Sciences) for constant interest in this work and stimulating discussions. He thanks Prof. N. A. Chernikov (Joint Institute for Nuclear Research, Dubna, Russia) for the interest in the problems posed in this work.

This research was partially supported by the Fund for Scientific Research of Bulgaria under contract Grant No. $F 103$.

\medskip
\medskip
 {\bf REFERENCES}

\medskip
1. Kobayashi S., K. Nomizu, Foundations of differential geometry, vol.1, Interscience publishers, New-York-London, 1963.\par
2. Schouten J. A., Ricci-Calculus: An Introduction to Tensor Analysis and its Geometrical Applications, Springer Verlag, Berlin-G\"otingen-Heidelberg, $2-nd ed., 1954$.\par
3. Iliev B. Z., Parallel transports in tensor spaces generated by derivations of tensor algebras, Communication
JINR, $E5-93-1$, Dubna, 1993.\par
4. Iliev B. Z., The equivalence principle in spaces with a linear transport along paths., Proceedings of the 5-th seminar "Gravitational energy and gravitational waves", JINR, Dubna, $18-20$ may 1992, Dubna, 1993, (In Russian).\par
5. Viro O.Ya., D.B. Fuks, I. Introduction to homotopy theory, In: Reviews of science and technic, sec. Modern problems in mathematics. Fundamental directions, vol.24, Topology-2, VINITI, Moscow, $1988, 6-121 ($In Russian).\par
6. Greub W., S. Halperin, R. Vanstone, Connections, Curvature, and Cohomology, vol.1, Academic Press, New York and London, 1972.\par
7. Hartman Ph., Ordinary differential equations, John Wiley\&Sons, New-York-London-Sydney, 1964.\par
8. Bellman R., Introduction to matrix analysis, McGRAW-HILL book comp., New York-Toronto-London, 1960.

\newpage
\vspace{5ex}
\noindent
 Iliev B. Z.\\[5ex]

\noindent
Linear Transports along Paths in Vector Bundles\\
 I. General Theory\\[5ex]

\medskip
\medskip
The (parallel) linear transports along paths in vector bundles are
axiomatically described. Their general form and certain properties are found.
It is shown that these transports are locally (i.e. along every fixed path)
always Euclidean ones in a senses that there exist frames in which their
matrices are unit. The investigated transports along paths are described in
terms of their local coefficients, as well as in terms of derivations along
paths.\\[5ex]

\medskip
\medskip
The investigation has been performed at the Laboratory of Theoretical
Physics, JINR.

\end{document}